\documentclass[11pt]{article}

\usepackage{amsmath}
\usepackage{amssymb}
\usepackage{eucal}

\begin{document}

\baselineskip 16pt

\title{On $\Pi$-permutable subgroups of   finite groups\thanks{Research is supported by
a NNSF grant of China (Grant \# 11371335) and Wu Wen-Tsun Key Laboratory of Mathematics of Chinese Academy of Sciences.} }

\author{Wenbin Guo\\
{\small Department of Mathematics, University of Science and
Technology of China,}\\ {\small Hefei 230026, P. R. China}\\
{\small E-mail:
wbguo@ustc.edu.cn}\\ \\
Alexander N. Skiba\\
{\small Department of Mathematics,  Francisk Skorina Gomel State University,}\\
{\small Gomel 246019, Belarus}\\
{\small E-mail: alexander.skiba49@gmail.com}}

\date{}
\maketitle

\begin{abstract} Let   $\sigma =\{\sigma_{i} | i\in I\}$ be   some
partition of the set of all primes $\Bbb{P}$ and    $\Pi$   a
non-empty   subset of the set $\sigma$.
A set  ${\cal H}$ of subgroups of a finite group $G$ is said to be a  \emph{ complete
Hall $\Pi
$-set} of $G$  if every  member of  ${\cal H}$ is a Hall $\sigma _{i}$-subgroup of $G$ for some $\sigma _{i}\in \Pi$ and    ${\cal H}$ contains exact one Hall  $\sigma
_{i}$-subgroup of $G$ for every  $\sigma _{i}\in \Pi$ such that $\sigma_i\cap \pi(G)\neq\emptyset$.
A subgroup $H$ of $G$ is called  \emph{$\Pi$-quasinormal} or
 \emph{$\Pi$-permutable} in $G$ if $G$ possesses a complete Hall $\Pi$-set
 ${\cal H}=\{H_{1}, \ldots ,
H_{t} \}$ such that $AH_{i}^{x}=H_{i}^{x}A$ for any $i$ and all $x\in G$.
We study the embedding  properties of $H$ under the hypothesis  that $H$ is $\Pi$-permutable in $G$.
Some known results are generalized.

\end{abstract}


\let\thefootnote\thefootnoteorig

\section{Introduction}

Throughout this paper, all groups are finite and $G$ always denotes
a finite group.  Moreover,
 $\Bbb{P}$ is the set of all    primes,
  $\pi \subseteq  \Bbb{P}$ and  $\pi' =  \Bbb{P} \setminus \pi$. If $n$ is an integer, the symbol $\pi (n)$
denotes the
 set of all primes dividing $|n|$; as usual,  $\pi (G)=\pi (|G|)$, the set of all
 primes dividing the order of $G$.

In what follows, $\sigma =\{\sigma_{i} | i\in I\}$ is  some
partition of $\Bbb{P}$, that is,
$\Bbb{P}=\cup_{i\in I} \sigma_{i}$ and $\sigma_{i}\cap
\sigma_{j}= \emptyset  $ for all $i\ne j$;   $\Pi$ is always supposed to be a
non-empty   subset of the set $\sigma$  and $\Pi'=\sigma\setminus
\Pi$.

In practice, we often deal with two limited  cases:   $\sigma =\{\{2\}, \{3\}, \{5\}, \ldots \}$ and    $\sigma =\{\pi, \pi'\}$.

Recall that  $\sigma (G)
=\{\sigma_{i} |\sigma_{i}\cap \pi (G)\ne  \emptyset  \}$  \cite{1}.   $G$ is called: a \emph{$\Pi$-group} if $\sigma (G) \subseteq \Pi$;
\emph{$\sigma$-primary} \cite{2} if $G$ is a $\Pi$-group for some one-element set $\Pi$.

A set  ${\cal H}$ of subgroups of $G$ is said to be a  \emph{ complete
Hall $\Pi
$-set} of $G$  if   every  member of  ${\cal H}$ is a Hall $\sigma _{i}$-subgroup of $G$ for some $\sigma _{i}\in \Pi$ and ${\cal H}$ contains exact one Hall  $\sigma
_{i}$-subgroup of $G$ for every  $\sigma _{i}\in \Pi \cap \sigma (G)$.
  We say also that $G$ is: \emph{ $\Pi$-full} if
 $G$ possesses a \emph{complete Hall $\Pi$-set};  a \emph{ $\Pi$-full group of Sylow type} if every subgroup of $G$ is a $D_{\sigma _{i}}$-group for all
 $\sigma _{i}\in \Pi$.

Let  $\cal L$ be some non-empty  set of subgroups of $G$ and $E$ a subgroup of $G$.
 Then a subgroup $A$ of $G$ is called
 \emph{$\cal L$-permutable}  if $AH=HA$ for all $H\in {\cal L}$;  \emph{${\cal L}^{E}$-permutable} if
 $AH^{x}=H^{x}A$ for all $H\in {\cal L}$ and all $x\in E$.

If  $\cal S$ is a complete Sylow  $\pi$-set of $G$ (that is, every member of $\cal S$ is a Sylow $p$-subgroup for some $p\in \pi$ and    $\cal S$ contains
exact one  Sylow  $p$-subgroup  for every $p\in \pi$), then  an ${\cal L}^{G}$-permutable subgroup
 is called  \emph{$\pi$-permutable} or \emph{$\pi$-quasinormal } (Kegel \cite{3}) in $G$.
 The $\pi (G)$-permutable subgroups are  also called  \emph{$S$-permutable}
or \emph{$S$-quasinormal. }

In this note we study the following generalization of $\pi$-permutability.

{\bf Definition 1.1.} We say that a subgroup $H$ of $G$ is \emph{$\Pi$-quasinormal} or
 \emph{$\Pi$-permutable} in $G$ if $G$ possesses a complete Hall $\Pi$-set
${\cal H}$ such that $H$ is ${\cal H}^{G}$-permutable.

Before continuing, consider some examples.

{\bf Example 1.2.}  (1)  $G$ is called \emph{$\sigma$-soluble} \cite{2} if
every chief factor of $G$ is $\sigma$-primary. In view of Theorem A in
\cite{1}, every  $\sigma$-soluble group is a  $\Pi$-full group of Sylow
type for each $\Pi \subseteq \sigma$.
.

(2)  $G$ is called \emph{$\sigma$-nilpotent} \cite{33} if  $G$ possesses a complete Hall $\sigma$-set ${\cal H}=\{H_{1}, \ldots ,
H_{t} \}$    such that  $G=H_{1}\times \cdots \times H_{t}$. Therefore every subgroup of every  $\sigma$-nilpotent group $G$ is $\Pi$-permutable in $G$ for each  $\Pi \subseteq \sigma$.

(3) Now let $p > q > r$ be   primes, where  $q$ divides $p-1$ and $r$ divides $q-1$.
 Let   $H=Q\rtimes R$
 be a non-abelian group of order $qr$, $P$  a simple    ${\mathbb F}_{p}H$-module  which is     faithful  for $H$,
 and  $G=P\rtimes H$.
Let $\sigma =\{\sigma _{1}, \sigma _{2}  \}$,  where  $\sigma _{1}= \{p, r\}$ and
$\sigma _{2} =\{p, r\}'$.
 Then $G$ is not $\sigma$-nilpotent and $|P| > p$. Since $q$ divides
$p-1$, $PQ$ is supersoluble. Hence for some normal subgroup  $L$ of
$PQ$ we have $1 < L < P$.
Then  for every Hall  $\sigma _{1}$-subgroup
$V$ of $G$ we have $L\leq P\leq V$, so $LV=V=VL$. On the other hand, for every
 Hall  $\sigma _{2}$-subgroup
$Q^{x}$ of $G$ we have $Q^{x}\leq PQ$, so $LQ^{x}=Q^{x}L$. Hence
$L$ is   $\sigma$-permutable in $G$. It is also  clear
 that $L$ is not normal in $G$, so
$LR\ne RL$, which implies that $L$  is   not
$S$-permutable in  $G$.

We will also need the following modification of the main concept in \cite{5}: A  subgroup $A$ of $G$ is called:   \emph{${\sigma}$-subnormal}
  in $G$ \cite{2}  if
there is a subgroup chain  $$A=A_{0} \leq A_{1} \leq \cdots \leq
A_{n}=G$$  such that  either $A_{i-1}$   is normal in $A_{i}$ or
$A_{i}/(A_{i-1})_{A_{i}}$ is  ${\sigma}$- primary
  for all $i=1, \ldots t$.

In this definition $(A_{i-1})_{A_{i}}$ denotes the product of all normal subgroups of $A_{i}$
contained in $A_{i-1}$.

We use $G^{{N}_{\sigma}}$ to denote  the \emph{$\sigma$-nilpotent residual} of $G$, that is,
 the intersection of all normal
subgroups $N$ of $G$ with $\sigma$-nilpotent quotient $G/N$.

Our main goal here is to prove the following

{\bf Theorem 1.3.}    {\sl  Let  $H$ be a $\Pi$-subgroup of $G$
 and    $D=G^{{N}_{\sigma}}$.}

(i) {\sl If $G$ is $\Pi$-full and possesses a complete Hall $\Pi$-set
  ${\cal H}$ such that $H$ is  ${\cal H}^{D}$-permutable, then
$H$ is $\sigma$-subnormal in $G$  and the normal closure    $H^{G}$ of $H$ in $G$  is a   $\Pi$-group.}

(ii)  {\sl  If $H$ is  $\Pi$-permutable in $G$  and, in the case when $\Pi\ne \sigma (G)$, $G$    possesses a complete Hall $\Pi'$-set   ${\cal K}$  such that $H$ is    ${\cal K}$-permutable,   then
$H^{G}/H_{G}$ is $\sigma$-nilpotent and the normalizer $N_{G}(H)$ of $H$
is also $\Pi$-permutable.  Moreover,  $N_{G}(H)$ is ${\cal H}^{G}$-permutable for each complete Hall $\Pi$-set
  ${\cal H}$ of $G$ such that $H$ is  ${\cal H}^{G}$-permutable.}

(iii)   {\sl If $G$ is a $\Pi'$-full group of  Sylow type  and
 $H$ is  $\Pi'$-permutable in $G$,
 then  $H^{G}$ possesses a $\sigma$-nilpotent Hall    $\Pi'$-subgroup.    }

Consider some corollaries of Theorem 1.3.

  Theorem 1.3(i)  immediately implies

 {\bf Corollary 1.4} (Kegel  \cite{5}).  {\sl    If  a $\pi$-subgroup
$H$ of $G$  is  $S$-permutable in $G$, then
$H$ is subnormal in $G$.
}

Now, consider some special cases of   Theorem 1.3(ii). First note that in
 the case when $\sigma =\{\{2\}, \{3\}, \ldots  \}$ we get from Theorem 1.3(ii) the following results.

 {\bf Corollary 1.5.} {\sl Let $H$ be a $\pi$-subgroup of $G$. If $H$ is $\pi$-permutable in $G$ and, also, $H$ permutes with some Sylow $p$-subgroup of $G$ for each prime $p\in \pi'$, then the normalizer $N_{G}(H)$ of $H$ is $\pi$-permutable in $G$.}

In particular, in the case when $\pi =\Bbb{P}$,  we have

{\bf Corollary 1.6} (Schmid \cite{99}).  {\sl If  a subgroup $H$ of $G$
is  $S$-permutable in $G$, then the normalizer
$N_{G}(H)$ of $H$  is also  $S$-permutable.}

{\bf Corollary 1.7.} {\sl Let $H$ be a $\pi$-subgroup of $G$. If $H$ is $\pi$-permutable in $G$ and, also, $H$ permutes with some Sylow $p$-subgroup of $G$ for each prime $p\in \pi'$, then $H/H_{G}$ is nilpotent.}

{\bf Corollary 1.8} (Deskins \cite{6}).  {\sl If a subgroup $H$
 of $G$ is $S$-permutable in $G$, then $H/H_{G}$ is nilpotent.}

Recall that  $G$ is said to be  a \emph{$\pi$-decomposable} if $G=O_{\pi}(G)\times O_{\pi'}(G)$,
 that is, $G$ is the direct
product of its  Hall $\pi$-subgroup and Hall  $\pi'$-subgroup.

In the case when $\sigma =\{\pi, \pi'\}$ we get from Theorem 1.3(ii) the following

{\bf Corollary 1.9.} {\sl Suppose that $G$ is $\pi$-separable. If a subgroup $H$ of $G$ permutes with all
Hall $\pi$-subgroups of $G$  and with  Hall $\pi'$-subgroups of $G$, then $H^{G}/H_{G}$ is $\pi$-decomposable.  }

In particular, we have

{\bf Corollary 1.10.} {\sl Suppose that $G$ is $p$-soluble. If a subgroup $H$ of $G$ permutes with all
Sylow  $p$-subgroups of $G$  and with   all  $p$-complements of $G$, then $H^{G}/H_{G}$ is $p$-decomposable.  }

Finally, in the case when $\Pi =\sigma$, we get from  Theorem 1.3(ii) the following

{\bf Corollary 1.11} (Skiba \cite{2}).  {\sl Suppose that $G$ is a $\sigma$-full group
 and let $H$ be a     subgroup of $G$.   If $H$ is  $\sigma$-permutable in $G$, then
$H^{G}/H_{G}$ is $\sigma$-nilpotent.}

 From Theorem 1.3(iii) we  get

{\bf Corollary 1.12.}  {\sl Let $H$ be a $\pi$-subgroup of $G$. If
$H$ permutes with every Sylow  $p$-subgroup of $G$ for $p\in \pi'$,  then  $H^{G}$  possesses a nilpotent
 $\pi$-complement.}

A subgroup $H$ of $G$ is called a \emph{$S$-semipermutable} in $G$ if $H$
permutes with all Sylow subgroups $P$ of $G$ such that $(|H|, |P|)=1$.
If    $H$ is $S$-semipermutable in $G$ and $\pi =\pi (H)$, then $H$ is $\pi'$-permutable in $G$.
Hence from Corollary 1.12 we  get the following known result.

{\bf Corollary 1.13} (Isaacs \cite{7}). {\sl If
  a  $\pi$-subgroup  $H$ of $G$ is  $S$-semipermutable in $G$,
 then  $H^{G}$  possesses a nilpotent
 $\pi$-complement.}

Note  that in the group $G=C_{7}\rtimes \text{Aut} (C_{7})$  a subgroup of order 3  is  $\pi'$-permutable  in $G$, where $\pi= \{2, 3\} $, but it is not $S$-semipermutable.

\section{Preliminaries}

We use:   $O^{\Pi}(G)$  to denote the subgroup of $G$ generated by all
its    ${\Pi}'$-subgroups;   $O_{\Pi}(G)$  to denote the subgroup of $G$ generated by all
its  normal   ${\Pi}$-subgroups.
 A  subgroup   $H$  of $G$ is said to be: a \emph{  Hall
 $\Pi$-subgroup} of  $G$ [1]  if $|H|$ is
a \emph{ $\Pi$-number}  (that is, $\pi(H)\subseteq \bigcup_{\sigma_i\in\Pi}\sigma_i$)  and $|G:H|$ is a  $\Pi'$-number.

{\bf Lemma 2.1.} {\sl Let  $A$,  $K$ and $N$ be subgroups of  $G$.
 Suppose that   $A$
is $\sigma$-subnormal in $G$ and $N$ is normal in $G$.  }

(1) {\sl $A\cap K$    is  $\sigma$-subnormal in   $K$}.

(2) {\sl If  $K$ is a $\sigma$-subnormal subgroup of  $A$,
then $K$ is $\sigma$-subnormal in $G$}.

(3) {\sl If $K$ is  $\sigma$-subnormal in $G$, then $A\cap K$  and  $\langle A, K\rangle$  are     $\sigma$-subnormal in
$G$.}

(4) {\sl $AN/N$ is $\sigma$-subnormal in $G/N$. }

(5) {\sl If $N\leq K$ and $K/N$ is $\sigma$-subnormal in $G/N$, then $K$
is $\sigma$-subnormal in $G$}.

(6) {\sl   If $K\leq A$ and $A$ is $\sigma$-nilpotent, then $K$
is $\sigma$-subnormal in $G$.}

(7) {\sl If $H\ne 1 $ is a Hall $\Pi$-subgroup of $G$  and $A$ is not  a  $\Pi'$-group, then $A\cap H\ne 1$ is
 a Hall $\Pi$-subgroup of $A$. }

(8) {\sl If  $|G:A|$ is a $\Pi$-number,  then  $O^{{\Pi}}(A)= O^{{\Pi}}(G)$.}

(9)    {\sl  If $G$ is   $\Pi$-full   and $A$ is a  $\Pi$-group,   then $ A\leq  O_{\Pi}(G)$. }

{\bf Proof.}  Statements (1)--(8) are known  \cite[Lemma 2.6]{2}).

(9)     Assume that this assertion   is false and let $G$ be a counterexample of
 minimal     order.   By hypothesis, there is a subgroup chain  $A=A_{0} \leq
A_{1} \leq \cdots \leq A_{r}=G$ such that
either $A_{i-1}$   is normal in $A_{i}$
  or $A_{i}/(A_{i-1})_{A_{i}}$ is  $\sigma $-primary  for all $i=1, \ldots , r$.
  Let   $M=A_{r-1}$.
  We can assume without loss of generality that $M\ne G$.
Let  $D=A\cap M_{G}$.

First   note that $A$ is not $\sigma$-primary. Indeed, assume that $A$ is a
$\sigma _{i}$-group.    By
hypothesis, $G$ has a Hall ${\sigma _{i}}$-subgroup, say $H$.  Then, by
Assertion (7), for any $x\in G$ we have $A\leq H^{x}$. Hence $A^{G}\leq
H_{G}\leq O_{\Pi}(G)$, a contradiction. Hence $|\sigma (A)| > 1$.

Suppose that  $D\ne 1$.  The subgroup  $D$ is $\sigma$-subnormal in $M_{G}$ by Lemma 2.1(1)(3),
 so the choice of $G$ implies that $D\leq O_{\Pi}(M_{G})$. Hence $O_{\Pi}(M_{G})\ne 1$.
But since  $O_{\Pi}(M_{G})$ is  characteristic in $M_{G}$, we have that $O_{\Pi}(M_{G})
\leq  O_{\Pi}(G)$. The hypothesis holds for $(G/O_{\Pi}(G), AO_{\Pi}(G)/O_{\Pi}(G))$
 by Assertion (4). Therefore  $AO_{\Pi}(G)/O_{\Pi}(G))\leq
O_{\Pi}(G/O_{\Pi}(G))=1$. It follows that $A\leq O_{\Pi}(G)$, a contradiction.
 Hence  $A\cap M_{G}=1$, so $M$ is not normal in $G$.
Therefore, $G/M_{G}$ is a $\sigma _{j}$-group for some  $j\in I$.
But then $A\simeq AM_{G}/M_{G}$ is $\sigma$-primary. This contradiction
completes the proof.

The first three statements in the next lemma can be proved by the direct calculations and the last statement see \cite[A, 1.6(a)]{DH}.

{\bf Lemma 2.2.}   {\sl  Let  $H$,  $K$ and $N$ be subgroups of $G$.
Let ${\cal H}=\{H_{1}, \ldots ,
H_{t} \}$ be a complete Hall $\Pi$-set of
 of $G$ and ${\cal L}={\cal H}^{K}$.
Suppose that   $H$
is ${\cal L}$-permutable and $N$ is normal in $G$. }

(1)  {\sl If  $H\leq E \leq G$, then $H$ is ${\cal L}^{*}$-permutable, where
  ${\cal L}^{*}=\{H_{1}\cap E, \ldots ,
H_{t}\cap E \}^{K\cap E}$. In particular, if  $H$ is
$\Pi$-permutable in $G$ and either  $G$ is a $\Pi$-full  group of
 Sylow type or $E$ is normal in $G$, then $H$ is $\Pi$-permutable in $E$.}

(2)  {\sl The subgroup  $HN/N$ is  ${\cal L}^{**}$-permutable,
 where ${\cal L}^{**}=\{H_{1}N/N, \ldots ,
H_{t}N/N \}^{KN/N}$.}

(3) {\sl If    $G$ is  a $\Pi$-full group of Sylow type and $E/N$ is a $\Pi$-permutable
subgroup of  $G/N$, then $E$
is $\Pi$-permutable in $G$.}

(4)  {\sl If   $K$ is $\cal L$-permutable, then  $\langle H, K\rangle$
 is $\cal L$-permutable.}

{\bf Lemma  2.3} (See Lemma 2.2 in \cite{1}). {\sl Let  $H$ be a normal subgroup of $G$. If
 $H/H\cap \Phi (G)$ is a $\Pi$-group, then $H$ has a  Hall $\Pi$-subgroup,
 say  $E$, and $E$ is normal in $G$}.

We say that a group $G$ is $\Pi$-closed if $O_{\Pi}(G)$ is a Hall $\Pi$-subgroup of $G$. Two integers $n$ and $m$ are called \emph{$\sigma$-coprime} if $\sigma (n) \cap \sigma (m)= \emptyset$.

{\bf Lemma 2.4.}    {\sl  If a $\sigma$-soluble  group $G$ has
   three  $\Pi$-closed  subgroups $A$, $B$ and $C$ whose indices
$|G:A|$, $|G:B|$, $|G:C|$ are pairwise  $\sigma$-coprime, then
$G$ is $\Pi$-closed. }

 {\bf Proof.}  Suppose  that this lemma   is false and let $G$ be a
counterexample with $|G|$ minimal. Let  $N$ be
 a minimal normal subgroup of $G$. Then  the hypothesis holds
 for $G/N$, so $G/N$ is $\Pi$-closed  by the choice of $G$. Therefore $N$
is not a  $\Pi$-group. Moreover,  $N$ is the unique minimal normal subgroup of $G$
 and, by Lemma 2.3,  $N\nleq \Phi (G)$.
Hence $C_{G}(N)\leq N$.  Since $G$ is   $\sigma$-soluble by
hypothesis, $N$ is $\sigma$-primary, say $N$ is a $\sigma _{i}$-group.
Then  $\sigma _{i}\in \Pi'$.

Since $|G:A|$, $|G:B|$, $|G:C|$ are pairwise  $\sigma$-coprime, there are at
least two subgroups, say $A$ and $B$, such  that $N\leq A\cap B$.  Then
$O_{\Pi}(A)\leq C_{G}(N)\leq N$, so  $O_{\Pi}(A)=1$.  But by hypothesis,
$A$ is $\Pi$-closed,  hence  $A$ is a  $\Pi'$-group.
Similarly we get that  $B$ is a  $\Pi'$-group and so $G=AB$ is
 a  $\Pi'$-group. But then $G$ is   $\Pi$-closed. This
contradiction completes the proof of the lemma.

Recall that $G$ is called a \emph{Schmidt group} if  $G$ is not nilpotent but every proper subgroup of $G$ is nilpotent.

{\bf Proposition 2.5.} {\sl  Let $G$ be a $\sigma$-soluble group.
Suppose that $G$    is not $\sigma  _{i}'$-closed  but all proper
subgroups of $G$  are $\sigma  _{i}'$-closed.  Then $G$ is  a  $\sigma  _{i}$-closed
  Schmidt group.    }

{\bf Proof.}   Suppose that this proposition is false and let $G$ be a
counterexample of minimal order. Let  $R$ be
 a minimal normal subgroup of   $G$ and $\{H_{1}, \ldots ,
H_{t} \}$    a  complete Hall $\sigma$-set of $G$.
Without loss of generality we can assume that $H_{1}$ is a  $\sigma
_{i}$-group.

(1)  {\sl   $|\sigma (G)|=2$. Hence $G=H_{1}H_{2}$.}

It is clear that $|\sigma (G)| > 1$. Suppose that $|\sigma (G)| > 2 $.
Then, since $G$ is $\sigma$-soluble, there are maximal subgroups $M_{1}$,
$M_{2}$ and $M_{3}$ whose indices $|G:M_{1}|$, $|G:M_{2}|$  and
$|G:M_{3}|$   are $\sigma$-coprime.
Hence $G=M_{1}M_{2}=M_{2}M_{3}=M_{1}M_{3}$. But the  subgroups $M_{1}$,
$M_{2}$ and $M_{3}$ are  $\sigma  _{i}'$-closed   by hypothesis. Hence
 $G$ is $\sigma  _{i}'$-closed  by Lemma 2.4, a contradiction.
Thus $|\sigma (G)|=2$.

(2) {\sl  If either $R\leq \Phi (G)$ or $R\leq H_{2}$, then  $G/R$ is a  $\sigma  _{i}$-closed
Schmidt group. }

 Lemma 2.3 and the choice of $G$ imply that  $G/R$ is not $\sigma  _{i}'$-closed.    On the other hand, every
maximal subgroup $M/R$ of $G/R$ is   $\sigma _{i}'$-closed since $M$ is $\sigma
_{i}'$-closed. Hence the hypothesis holds for $G/R$. The
choice of $G$   implies that   $G/R$ is a  $\sigma  _{i}$-closed
Schmidt group.

(3)  {\sl  $\Phi (G)=1 $,  $R$ is  the unique
minimal normal subgroup of $G$ and $R\leq  H_{1}$.   }

Suppose that    $R\leq \Phi (G)$.  Then $R$ is a $r$-group for some prime
$r$ and,  in view of Claim (1),  Lemma 2.3  and \cite[IV, 5.4]{hupp},  $G=
H_{1}\rtimes H_{2}=P\rtimes Q$, where  $H_{1}=P$ is a $p$-group and $H_{2}=Q$
 is a $q$-group for some different primes $p$ and $q$.  Assume that $R\leq Q$ and take a subgroup $L$ of order $q$ in
$R\cap  Z(Q)$. Then it is clear that $R < Q$, so $PR < G$  and hence
$PR=P\times Q$ is $p$-nilpotent. Therefore $L\leq  Z(G)$, so  $R=L\leq Z(G)$. But
 for every maximal subgroup  $M$ of $G$ we have
 $R\leq M$ and  $M/R$ is nilpotent. Hence every maximal subgroup of $G$ is
nilpotent and  so $G$ is a $\sigma  _{i}$-closed
Schmidt group, a contradiction.  Similarly, we get that $G$ is a $\sigma  _{i}$-closed
Schmidt group in the case when $R\leq P$.  Therefore $R\nleq \Phi (G)$.

Now  assume that $G$ has a minimal normal subgroup $L\ne R$.
Then by (3), there are maximal subgroups $M$ and $T$   of $G$ such that $LM=G$ and
$RT=G$. By hypothesis, $M$ and $T$ are  $\sigma  _{i}'$-closed. Hence
$G/L\simeq LM/L\simeq M/M\cap L$ is $\sigma  _{i}'$-closed. Similarly,
$G/R$ is $\sigma  _{i}'$-closed and so $G\simeq G/L\cap R$ is $\sigma
_{i}$-nilpotent, a contradiction.  Hence $R$ is  the unique
minimal normal subgroup of $G$,  and so $R\leq  H_{1}$.

  {\sl  Final contradiction.}  In view of Claim (3),
$C_{G}(R) \leq R$. Hence $|H_{2}|$ is a prime and
   $RH_{2}=G$ since $R\leq H_{1}$ and every proper subgroup of $G$
 is $\sigma  _{i}'$-closed. Therefore $R=H_{1}$, so $R$ is not abelian since  $G$ is a not
a $\sigma  _{i}$-closed  Schmidt group.    By Claim (1) and Theorem 3.5 in \cite{Gor}, for any prime
$p$ dividing $|R|$ there is a Sylow $p$-subgroup $P$ of $G$ such that $PH_{2}=H_{2}P$.
But $H_{2}P < G$, so   $H_{2}P=H_{2}\rtimes  P$. This implies that $R\leq
N_{G}(H_{2})$ and thereby $G=R\times H_{2} =H_{1}\times H_{2}
$.  This final contradiction completes the proof of the result.

{\bf Corollary  2.6.}    {\sl Let $G$ be a minimal non-$\sigma$-nilpotent group,  that is, $G$ is not $\sigma$-nilpotent, but every proper subgroup of $G$ is $\sigma$-nilpotent. If $G$ is a $\sigma$-soluble, then $G$ is a Schmidt group.}

{\bf Proof.} It is clear that $G$ is $\sigma$-nilpotent if and only if $G$ is $\sigma _{i}'$-closed for all $\sigma _{i}\in \sigma$. Hence, for some $i$, $G$ is not $\sigma _{i}'$-closed. On the other hand, every proper subgroup of $G$ is $\sigma _{i}'$-closed. Hence  $G$ is a Schmidt group by Proposition 2.5.

{\bf Proposition  2.7.}    {\sl  Let $G$ be a $\Pi$-full  group of Sylow type. If $G$
   possesses a $\sigma$-nilpotent Hall
 $\Pi$-subgroup $H$, then every $\sigma$-soluble
  $\Pi$-subgroup  of $G$ is contained in a conjugate of $H$.
 In particular, any two $\sigma$-soluble Hall
 $\Pi$-subgroups of $G$ are conjugate.   }

{\bf Proof.} Suppose that this proposition is false  and let  $G$ be a
counterexample of minimal order.  Then  some $\sigma$-soluble $\Pi$-subgroup $K$ of $G$ is
not contained in $H^{x}$ for all $x\in G$. We   can assume without loss of
generality that every proper subgroup $V$ of $K$ is contained in a conjugate of $H$, so
$V$ is $\sigma$-nilpotent. Hence either $K$ is $\sigma$-nilpotent or $K$ is a minimal non-$\sigma$-nilpotent group.
Then in view of Corollary 2.6   and \cite[IV, 5.4]{hupp},  $K$ has a
 normal Hall $\sigma _{i}$-subgroup $L$ for some $\sigma _{i}\in \sigma (K)$.
Now arguing as in the proof of Wielandt's theorem \cite[(10.1.9)]{Rob},
one can show that for some $y\in G$ we have $K\leq H^{y}$. This
contradiction completes the proof of the result.

{\bf Corollary 2.8.}    {\sl   Let $G$ be a $\Pi$-full  group of Sylow type. Suppose that every chief factor of $G$ possesses  a $\sigma$-nilpotent Hall $\Pi$-subgroup. Then $G$ possesses a $\sigma$-soluble Hall $\Pi$-subgroup.
   }

{\bf Proof.} Let $R$ be a minimal normal subgroup of $G$,   $H$ a $\sigma$-nilpotent Hall $\Pi$-subgroup of $R$ and $N=N_{G}(H)$.
 By induction, $G/R$ has a  $\sigma$-soluble  Hall $\Pi$-subgroup, say $U/R$. Therefore if $R$ is a $\Pi$-group, then $U$ is a $\sigma$-soluble Hall $\Pi$-subgroup of $G$. On the other hand, if $R$ is a $\Pi'$-group, then $U=R\rtimes V$ by the Schur-Zassenhas theorem, where $V\simeq U/R$ is a $\sigma$-soluble Hall $\Pi$-subgroup of $G$.  Now suppose that
$1 < H < R$.  Proposition 2.7 and the Ftattini argument imply that $G=RN$, where   $|G:N|=|R/R\cap N|$ is a $\Pi'$-number and $N < G$.
Then   $N/N\cap R\simeq G/R$  possesses a $\sigma$-soluble Hall $\Pi$-subgroup. Hence in view Proposition 2.7, the hypothesis holds for $N$, so $N$ possesses a $\sigma$-soluble  Hall $\Pi$-subgroup $W$ by induction. It is clear now that $W$ is a Hall $\Pi$-subgroup of $G$.
The corollary is proved.

\section{Proof of Theorem 1.3}

Suppose  that this theorem is false and let $(G, H)$ be a counterexample
with $|G| + |G:H|$ as small as possible. Then $H\ne H^{G}$.

(i), (ii)   By hypothesis, $G$ possesses  a complete Hall $\Pi$-set, say ${\cal H}= \{H_{1}, \ldots ,H_{t}\}$.  We can assume without loss
 of generality that $H_{i}$ is
a  $\sigma _{i}$-group for all $i=1,  \ldots , t$. Let    $E=H_{1}^{G}\cdots H_{t}^{G}$.

Suppose that Assertion (i) is false. Then in view of Lemma 2.1(9),  $H$ is not $\sigma$-subnormal in $G$.
Moreover, in this case we have  $E= G$. Indeed, since the class of all $\sigma$-nilpotent groups is closed under taking subgroups, homomorphic images and the direct products,   $E/E\cap D\simeq DE/D$ is $\sigma$-nilpotent. Hence $E^{{N}_{\sigma}}\leq D$. It follows that
 the hypothesis holds for $(E, H)$. Thus in the case when  $E < G$ the choice of $(G, H)$ implies that $H$ is $\sigma$-subnormal in $E$ and so $H$ is $\sigma$-subnormal in $G$, a contradiction. Therefore $E= G$. Since $H\ne H^{G}$, it follows that for some $x\in G$ and $H_{i}\in {\cal H}$
we have $H_{i}^{x}\nleq N_{G}(H)$. Now, arguing as in Claim (2) of the proof of Theorem B in \cite{2}, one can show that $H$ is $\sigma$-subnormal in $G$. This contradiction completes the proof of (i).

(ii) Suppose that this assertion is false. 
Then:

(1) {\sl  The hypothesis holds for $(G/H_{G}, H/H_{G})$, so $H_{G}=1$.}

First note that    the hypothesis holds for $(G/H_{G}, H/H_{G})$ by Lemma  2.2(2).
 Assume that $H_{G}\ne 1$. Then the choice of $(G, H)$ implies that
$H^{G}/H_{G}$ is $\sigma$-nilpotent and $N_{G/H_{G}}(H/H_{G})=  N_{G}(H)/H_{G}$ is ${\cal H}^{*}$-permutable by Lemma 2.2(2), where
$${\cal H}^{*}=\{H_{1}H_{G}/H_{G}, \ldots , H_{t}H_{G}/H_{G}\}^{G/H_{G}}.$$ But then,  clearly, $N_{G}(H)$ is ${\cal H}^{G}$-permutable.
This shows that Assertion (ii) is  true. Therefore the choice of $(G,H)$ implies that $H_{G}=1$.

(2)  {\sl $t > 1$.}

Assume that $t=1$, that is, $H$ is  a $\sigma _{1}$-group.
Then   $HH_{1}^{x}=H_{1}^{x}H=H_{1}^{x}$  for all $x\in G$,  so $H^{G}\leq (H_{1})_{G}\leq  O_{\sigma _{1}}(G)$, which implies that $H^{G}$ is $\sigma$-nilpotent. Hence  $H$ is
$\sigma$-subnormal in $G$ by Lemma 2.1(6). Note also that for any Hall $\sigma _{1}'$-subgroup $V$ of $G$ such that $HV=VH$ we have $H=VH\cap O_{\sigma _{1}}(G)$, so $V\leq N_{G}(H)$. Therefore  if $H$ is  $\Pi$-permutable in $G$  and also, in the case when $\Pi\ne \sigma (G)$,  $H$ is ${\cal K}$-permutable,   then $|G:N_{G}(H)|$ is a  $\sigma _{1}$-number, which implies that $N_{G}(H)H_{1}^{x}=G=H_{1}^{x}N_{G}(H)$ for all $x\in G$. This means that $N_{G}(H)$ is $\Pi$-permutable in $G$. Thus Assertion (ii) is true, a contradiction. Therefore $t > 1$.

 Let $L_{i}=O^{\sigma _{i}'}(H)$, for all $i=1, \ldots , t$.
 Then $H=L_{1}\cdots L_{t}$ and
$N_{G}(H)=N_{G}(L_{1})  \cap  \cdots   \cap N_{G}(L_{t})$. Let
 $$W_{i}=H_{1}^{G} \cdots H_{i-1}^{G}H_{i+1}^{G} \cdots H_{t}^{G},$$ for all $i=1, \ldots , t$, and $W=W_{1}\cap \cdots \cap W_{t}$.

(3) {\sl $W_{i} \leq N_{G}(L_{i})$   for all $i=1, \ldots , t$, so $W\leq  N_{G}(H)$}.

Indeed, since $H$ is $\sigma$-subnormal in $G$ by Part (i),
Lemma 2.1(8) implies that $H_{i}^{x}\leq  N_{G}(O^{\sigma _{i}}(H))$ for all $x\in G$.
 This means that $H_{i}^{G}\leq  N_{G}(O^{\sigma _{i}}(H))$. Hence
$H_{i}^{G}\leq  N_{G}(L_{j})$  for all $j\ne i$, so $W_{i} \leq N_{G}(L_{i})$   for all $i=1, \ldots , t$.

 (4) $H^{G}$ is $\sigma$-nilpotent.

Suppose that this is false.  Let  $K=H_{1}\cdots H_{t}W$. Then:

(a)  {\sl $K$ is a  subgroup of $G$, $H\leq K$ and $|K:W|$ is a $\Pi$-number.}

First note that $(H_{i}W/W)^{G/W}=H_{i}^{G}W/W$ and

$$WW_{i}\cap H_{i}^{G}W=W(W_{i}\cap H_{i}^{G}W)=W(W_{i}\cap H_{i}^{G}(W_{1}\cap \cdots \cap W_{t}))=$$$$=W(W_{i}\cap W_{1}\cap  \cdots \cap W_{i-1}\cap W_{i+1}\cap  \cdots \cap W_{t} \cap W_{i}H_{i}^{G})= W(W\cap E)=W.$$
Therefore $$E/W=(H_{1}W/W)^{G/W} \times \cdots \times (H_{t}W/W)^{G/W}.$$
This means that $[H_{i}W/W, H_{j}W/W]=1$, for all $i\ne j$.
Hence  $K=H_{1}\cdots H_{t}W=(H_{1}W) \cdots  (H_{t}W)$ is the product of pairwise permutable subgroups, which implies that  $K$ is
a subgroup  of $G$.  It is also clear that $K/W$ is a Hall $\Pi$-subgroup of $G/W$. Hence $|K:W|$ is a $\Pi$-number and   $WH/W\leq  K/W$ by  Lemma 2.1(4)(7), so we have (a).

(b) {\sl The hypothesis holds for $(K, H)$.}

Let   ${\cal K}= \{K_{1}, \ldots , K_{n}\}.$  Since  $|K:W|$ is a $\Pi$-number, $K_{i}\cap K$ is a Hall
  $\sigma _{i}$-subgroup  of $K$ and   hence ${\cal B}=\{K_{1}\cap K, \ldots , K_{n}\cap K \}$ is a complete Hall $\Pi'$-set of $K$.
 On the other hand, for  any  $K_{i}\in \cal K$ we have $HK_{i}\cap K=(K_{i}\cap K)H$ and so $H$ is $\cal B$-permutable.
  Finally, it is clear that $H$ is $\Pi$-permutable in $K$. Hence  the hypothesis holds for  $(K, H)$.

(c) $K < G$.

 Suppose that $K=G$. Then, since $|K:W|=|G:W|$ is a $\Pi$-number by Claim (4),
for every $K_{i}\in \cal K$ and every $x\in G$ we have $K_{i}^{x} \leq  W \leq  N_{G}(H)$ by Claim (3), so $K_{i}^{x}H=HK_{i}^{x}$. Therefore $H$ is $\sigma$-permutable in $G$ and so $H^G\simeq H^{G}/H_{G}$ is $\sigma$-nilpotent by Theorem B in \cite{2}, contrary to our assumption on $H$. Hence $K < G$.

(d) {\sl $|G:N_{G}(H)|$ is a $\Pi$-number} (Since  $H$ is a $\sigma$-subnormal $\Pi$-subgroup of $G$, this follows from Lemma 2.1(8)).

(e) Conclusion for (4).

Since $K < G$ by Claim (c), we have that   $H^{K}/H_{K}$ is $\sigma$-nilpotent. Because $|G:N_{G}(H)|$ is a $\Pi$-number by Claim (d),  $G=KN_{G}(H)$. Hence $H^{G}\simeq H/1=H^{G}/H_{G}=H^{K}/H_{K}$ is $\sigma$-nilpotent. This contradiction shows that $H^{G}$ is $\sigma$-nilpotent.

{\it Final contradiction for (ii).}

Since $H^G$ is $\sigma$-nilpotent by (4), $H$ is also $\sigma$-nilpotent. Hence $H$ possesses a complete Hall $\sigma$-set $\{V_{1}, \ldots ,
V_{t} \}$    such that  $H=V_{1}\times \cdots \times V_{t}$.
Without loss of generality we can assume that $V_{i}$ is a $\sigma _{i}$-group for all $i=1, \ldots, t$.
 Let $N=N_{G}(H)$ and $N_{i}=N_{G}(V_{i})$. Then $N=N_{1}\cap \cdots \cap N_{t}$. Moreover, it is clear that $L_{i}=V_{i}$ for all $i=1, \ldots, t$.
Hence $W_{i}\leq  N_{G}(V_{i})$ for all $i=1, \ldots, t$ by Claim (3). It is also clear  that $|G:N_{i}|$ is a $\sigma _{i}$-number, so $G= N_{i}H_{i}$.
Hence for any $x\in G$ and $H_{i}\in \cal H$ we have $$NH_{i}^{x}=(N_{1}\cap \cdots \cap N_{t})H_{i}^{x}=N_{i}H_{i}^{x}\cap N_{1}\cap
\cdots \cap N_{i-1}\cap N_{i+1} \cap \cdots \cap N_{t}=$$$$=G\cap N_{1}\cap \cdots \cap N_{i-1}\cap N_{i+1} \cap \cdots \cap N_{t}=N_{1}\cap \cdots \cap N_{i-1}\cap N_{i+1} \cap \cdots \cap N_{t}=H_{i}^{x}N$$ and so $N$ is ${\cal H}^{G}$-permutable.
Therefore Assertion (ii) is true.
This contradiction completes the proof of Assertion (ii).

(iii)  Let   ${\cal
L}=\{L_{1}, \ldots , L_{m} \}$ be a complete Hall $\Pi'$-set  of  $G$ such that $H$ is ${\cal L}^{G}$-permutable.
Let $E=H^{G}$ and $R$ a minimal normal subgroup of $G$. First note that $m > 1$, Indeed, if $m=1$, then $L_{1}\cap E$ is a $\sigma$-nilpotent Hall $\Pi'$-subgroup of $G$, which contradicts the choice of $(G, H)$.

(1) {\sl $ER/R$ possesses
a $\sigma$-nilpotent Hall  $\Pi'$-subgroup $U/R$. Therefore $R\leq E$. }

From Lemma 2.2(2)  and the choice of $G$ it
 follows that $(HR/R)^{G/R}=   ER/R$  possesses a $\sigma$-nilpotent Hall  $\Pi'$-subgroup, say $U/R$. Therefore, if
$R\nleq E$, then $E\simeq ER/R$  possesses
a $\sigma$-nilpotent Hall  $\Pi'$-subgroup, a contradiction.
 Hence we have (1).

(2) {\sl $O_{\Pi}(G)=1$.}

Assume that  $R\leq O_{\Pi}(G)$. Then, by the Schur-Zassenhaus theorem, $R$ has a complement $V$ in
$U$, so $V\simeq U/R$ is a $\sigma$-nilpotent  Hall $\Pi'$-subgroup of
$E$, a contradiction.   Hence we have (2).

(3) {\sl  $L_{i}^{G}\nleq  C_{G}(E)$ for all $i=1, \ldots , t$.}

Assume that $L_{i}^{G}\leq  C_{G}(E)$ and let $N$ be a minimal normal subgroup
 of $G$ contained in    $L_{i}^{G}$. Then $N\leq E$ and $E/N$ possesses
a $\sigma$-nilpotent Hall  $\Pi'$-subgroup, say $U/N$, by Claim (1). On the other hand,
$N\leq Z(U)$, so $U$ is   $\sigma$-nilpotent.  But  a Hall
$\Pi'$-subgroup of $U$ is a Hall  $\Pi'$-subgroup of $E$, a
contradiction. Hence we have (3).

(4) {\sl $R$ is the unique minimal normal subgroup of $G$.}

Suppose that $G$ has a minimal normal subgroup $N\ne R$.
Then $N\leq E$   and  $G/N$ possesses  a  $\sigma$-nilpotent Hall
$\Pi'$-subgroup  by Claim (1). Therefore  $(E/R)\times (E/N)$ possesses a  $\sigma$-nilpotent Hall
$\Pi'$-subgroup $V$.  But  $E\simeq K\leq  (E/R)\times (E/N)$ since $R\cap N=1$. Hence $E$ possesses a $\sigma$-nilpotent Hall $\Pi'$-subgroup.
Moreover, since $N\simeq  RN/R$ possesses a $\sigma$-nilpotent Hall $\Pi'$-subgroup, $E$ possesses a Hall $\Pi'$-subgroup $U$ by Corollary 2.8.
  But then, by Proposition 2.7, for some $x\in G$ we have $U\leq V^{x}$
 and so $U$  is  $\sigma$-nilpotent, contrary to the choice of $G$. Hence
we have (4).

{\sl Final contradiction for (iii).}

Let  $x, y\in G$ and  $A=H^{x}$.
Then  $$AL^{y}_{i}=(HL^{yx^{-1}}_{i})^{x}=(L^{yx^{-1}}_{i}H)^{x}=L^{y}_{i}A$$ by hypothesis.
  Let $L=A^{L_{i}}\cap L_{i}^{A}$.  Then $L$ is a subnormal subgroup
of $G$ by  \cite[7.2.5]{stounh}. Suppose that $L\ne 1$ and let
 $L_{0}$ be a minimal subnormal subgroup of $G$ contained in $L$. Then $V=L_{0}\cap L_{i}  $ is a Hall $\Pi'$-subgroup of $L_{0}$  since   $L\leq AL_{i}$. Moreover, in view of Claim (2), $V\ne 1$ (see, for example, \cite[Chapter 1, Lemma 5.35(5)]{Guo}).
  We now show that   $L_{i}\cap R$ is a non-identity    Hall
$\Pi'$-subgroup  of $R$. Indeed, if $L_{0}$ is abelian, then $L_0 \leq O_{\sigma _{i}}(G)$, where $\sigma_i=\pi(L_i),$ so
$R$ is a $\sigma _{i}$-group by Claim (4). On the other hand, if $L_{0}$ is non-abelian, $L_{0}^{G}$ is a minimal normal subgroup of $G$ and  so, by Claim (4),  $L_{i} \cap R$  is a non-identity    Hall
$\Pi'$-subgroup  of $R$.

Since $m > 1$,  Claim (2) implies that  there
is $j\ne i$ such that   for every  $x, y\in G  $ we have
 $(L^{y}_{j})^{H^{x}}\cap (H^{x})^{L^{y}_{j}}=1$ and so
 $$[L^{y}_{j}, H^{x}]\leq [(L^{y}_{j})^{H^{x}}, (H^{x})^{L^{y}_{j}}] =1.$$
 Therefore  $L_{j}^{G}\leq  C_{G}(E)$, contrary  Claim (3).  Hence Statement (iii) holds.

 The theorem is proved.

\end{document}